\documentclass[11pt,a4paper]{amsart}

\usepackage{xy}

\usepackage{amsfonts}
\usepackage{amsthm}
\usepackage[T1]{fontenc}
\usepackage[utf8]{inputenc}
\usepackage{dsfont}
\usepackage[english]{babel}

\xyoption{all}

\usepackage{color}
\usepackage{graphics}
\usepackage{amssymb}

 \theoremstyle{plain}
 
 \newtheorem{Prop}{Proposition}[section]
 \newtheorem{Lem}[Prop]{Lemma}
 \newtheorem{Cor}{Corollary}
 \newtheorem*{mainthm}{Main Theorem}

 \theoremstyle{definition}

\newenvironment{pf}{\begin{proof}}{\end{proof}}

\newcommand{\upcirc}[1]{\stackrel{#1}{\circ}}

\newcommand{\iso}{\cong}

\newcommand{\im}{\mathrm{Im}}

\newcommand{\rep}[4]{\left( \left\{#1_{#3} \right\}_{#3\in #2_0}, \left\{#1_{#4} \right\}_{#4\in #2_1}\right)}

\newcommand{\kwadracik}{\hfill\qed}

\newcommand{\modplus}[1]{\mathrm{mod}_+(#1)}
\newcommand{\modminus}[1]{\mathrm{mod}_-(#1)}
\newcommand{\modzero}[1]{\mathrm{mod}_0(#1)}
\newcommand{\cohplus}[1]{\mathrm{coh}_+(#1)}
\newcommand{\cohminus}[1]{\mathrm{coh}_-(#1)}

\newcommand{\id}{\mathds{1}}

\newcommand{\Tt}{\mathcal{T}}

\newcommand{\Ff}{\mathcal{F}}

\newcommand{\Ll}{\mathcal{L}}
\newcommand{\Oo}{\mathcal{O}}

\newcommand{\XX}{\mathbb{X}}
\newcommand{\NN}{\mathbb{N}}
\newcommand{\LL}{\mathbb{L}}

\newcommand{\ZZ}{\mathbb{Z}}

\newcommand{\pp}{\underline{p}}
\newcommand{\lala}{\underline{\la}}

\newcommand{\al}{\alpha}
\newcommand{\La}{\Lambda}
\newcommand{\la}{\lambda}

\newcommand{\tauX}{\tau_\XX}
\newcommand{\lra}{\longrightarrow}

\newcommand{\End}[2]{\mathrm{ End}_{#1}\left(#2\right)}
\newcommand{\EndX}[1]{\End{\XX}{#1}}

\newcommand{\Hom}[3]{\mathrm{ Hom}_{#1}\big(#2,#3\big)}
\newcommand{\HomX}[2]{\Hom{\XX}{#1}{#2}}
\newcommand{\HomL}[2]{\Hom{\La}{#1}{#2}}

\newcommand{\Ext}[4]{\mathrm{ Ext}^{#1}_{#2}\big(#3,#4\big)}
\newcommand{\ExtX}[2]{\Ext{1}{\XX}{#1}{#2}}
\newcommand{\ExtL}[2]{\Ext{1}{\La}{#1}{#2}}

\newcommand{\coh}[1]{\mathrm{coh}(#1)}
\newcommand{\vect}[1]{\mathrm{vect}(#1)}

\newcommand{\cohnull}[1]{\mathrm{coh}_0(#1)}
\renewcommand{\mod}[1]{\mathrm{mod}(#1)}

\newcommand{\rk}{\mathrm{rk}}

\newcommand{\vx}{{\vec{x}}}
\newcommand{\vy}{\vec{y}}
\newcommand{\vc}{\vec{c}}

\newcommand{\vz}{{\vec{z}}}
\newcommand{\vw}{{\vec{\omega}}}

\begin{document}

\baselineskip=17pt

\title[Exceptional modules]{Exceptional modules over wild canonical algebras}

\author{Dawid Edmund Kędzierski}
\author{Hagen Meltzer}
\date{}


\address{Institut of Mathematics, Szczecin University, $70-451$
  Szczecin, Poland}

\email{dawid.kedzierski@usz.edu.pl,
hagen.meltzer@usz.edu.pl}





\begin{abstract}

We show that ''almost all'' exceptional modules over  wild canonical algebra $\La(\pp,\lala)$ can be described by matrices having coefficients $\la_i-\la_j$, where $\la_i, \la_j$ are elements from the parameter sequence $\lala$.

The proof is based on Schofield induction for sheaves in the
associated categories of weighted projective lines \cite{Kedzierski:Meltzer:2013}
and an extended version of C. M. Ringel's proof for the $''0, \,1''$ matrix
property for exceptional representations for finite acyclic quivers \cite{Ringel:1998}.

\end{abstract}

\maketitle

\vspace{-0,5cm}
\section{Introduction}
Canonical algebras were introduced by C. M. Ringel \cite{Ringel:1984}. A  \emph{canonical algebra} $\La$ of quiver type
over a field $k$ is a quotient algebra of the path algebra  of the  quiver $Q$:
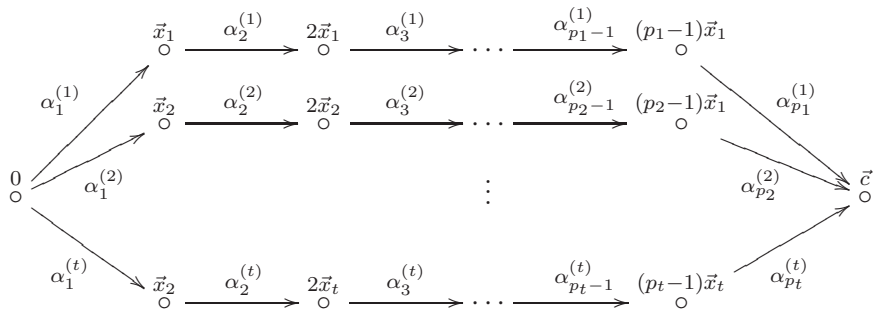
$$\xymatrix @C +1.5pc @R-1.5pc{&\upcirc{\vx_1}\ar[r]^{\al_2^{(1)}}& \upcirc{2\vx_1}\ar[r]^{\al_3^{(1)}}&\cdots \ar[r]^{\al_{p_1-1}^{(1)}}&\upcirc{(p_1-1)\vx_1} \ar[rdd]^{\al_{p_1}^{(1)}}&\\
&\upcirc{\vx_2}\ar[r]^{\al_2^{(2)}}& \upcirc{2\vx_2}\ar[r]^{\al_3^{(2)}}&\cdots \ar[r]^{\al_{p_2-1}^{(2)}}&\upcirc{(p_2-1)\vx_1}\ar[rd]_{\al_{p_2}^{(2)}}&\\
\upcirc{{0}}\ar[ruu]^{\al_1^{(1)}} \ar[ru]_{\al_1^{(2)}} \ar[rdd]_{\al_1^{(t)}}  &&&\vdots&&\upcirc\vc\\
&& && &\\
&\upcirc{\vx_{2}}\ar[r]^{\al_2^{(t)}}& \upcirc{2\vx_{t}}\ar[r]^{\al_3^{(t)}}&\cdots \ar[r]^{\al_{p_{t}-1}^{({t})}}&\upcirc{(p_{t}-1)\vx_{t}}\ar[ruu]_{\al_{p_{t}}^{({t})}}&\\
}$$
modulo the ideal $I$ defined by the \emph{canonical relations}
$$\al_{p_i}^{(i)}\dots\al_{2}^{(i)}\al_1^{(i)}=\al_{p_1}^{(1)}\dots \al_{2}^{(1)}\al_1^{(1)}+\lambda_i\al_{p_2}^{(2)}\dots \al_{2}^{(2)}\al_1^{(2)}\quad \text{for}\quad i=3,\dots,t,$$
where the $\la_i$ are pairwise distinct non-zero elements of $k$. They are called \emph{parameters}. The positive integer numbers $p_i$ are at least $2$ and they are called the  \emph{weights}.
Usually we assume that $k$ is algebraically closed, but for many results this is not necessary.
The algebra $\La$ depends on a weights sequence $\pp=(p_1,\dots, p_t)$ and a sequence of parameters $\lala= (\la_2,\dots,\la_t)$. We can assume that $\la_2=0$ and $\la_3=1$. We write $\La=\La(\pp,\lala)$. Concerning the complexity of the module category over $\La$ there are three types of canonical algebras: domestic, tubular and wild. Recall that the $\La$ is of wild type if and only if the Euler characteristic $\chi_\La= (2-t)+\sum_{i=1}^t1/{p_i}$ is negative.

Denote by $Q_0$ the set of vertices and by $Q_1$ the set of arrows of the quiver $Q$. Then each finite generated right module over $\La$ is given by finite-dimensional vector spaces $M_i$ for each vertex $i$ of $Q_0$ and by linear maps $M_\al:M_j\rightarrow M_i$ for the arrows $\al:i\rightarrow j$ of
$Q_1$ such that the canonical relations are satisfied.
 We will usually identify the linear maps with matrices. The category of all this modules we denote by $\mod\La$.

Our aim is to study the possible coefficients, which can appear in the matrices of exceptional modules over wild canonical algebras.
In many cases the matrices of special modules can be exhibited by $0$, $1$ matrices.
This was shown  by C. M. Ringel for exceptional representations of finite acyclic quivers
\cite{Ringel:1998}
and for indecomposable modules over representation-finite algebras, which is a result of P. Dr\"axler \cite{Draxler:2001}.
In some special case explicit $0$, $1$  matrices with few nonzero entries have been calculated,
so for indecomposable representation of Dynkin  quivers  by P. Gabriel \cite{Gabriel} and
indecomposable representation of representation-finite posets by M. Kleiner \cite{Kleiner}
(see also a result of K. J.  B\"ackstroem for orders over lattices \cite{Baeckstroem}).
Among new results we mention  a paper of M. Grzecza, S. Kasjan and A. Mróz \cite{Grzecza:Kasjan:Mroz:2012}.

The  problem of determining matrices for indecomposable modules over canonical algebras has  been solved in the case of domestic case.
In  the case of a field of characteristic
different from $2$ D. Kussin and the second author  computed  matrices having entries $0$, $\pm 1$ for all indecomposable modules, where the entries $-1$ appears only for very special regular modules \cite{Kussin:Meltzer:2007a}.
Matrices of indecomposable modules over canonical algebras over an arbitrary field were described in
 \cite{Komoda:Meltzer:2008}.
These results were used to determine matrices for exceptional representations for tame quivers
\cite{Kussin:Meltzer:2007b}, \cite{Kedzierski:Meltzer:2011}.

In the case of tubular canonical algebras it was shown in \cite{Meltzer:2007}  that each exceptional module can be described by matrices having as entries
$0$, $\pm 1$ in the tubular types $(2,3,6)$, $(3,3,3)$, $(2,4,4)$ and for the weight type $(2,2,2,2)$ with a parameter $\la$ appear entries $0$, $\pm 1$, $\pm \la$ and $1-\la$. The proof uses universal extensions in the sense of K. Bongartz \cite{Bongartz}.

Later P. Dowbor, A. Mróz and the second author  developed an algorithm and a computer program
for explicit calculations of matrices for  exceptional modules over tubular canonical algebras \cite{Dowbor:Meltzer:Mroz:2010}.

 In general  little is known about matrices of non-exceptional modules. However in the case of tubular canonical algebra
an algorithm for the computation of matrices of non-exceptional modules was developed in \cite{Dowbor:Meltzer:Mroz:2014bimodules}.
Moreover, explicit formulas for these matrices were obtained in case the module is of integral slope  \cite{Dowbor:Meltzer:Mroz:2014slope}.

Recently the $0$, $1$ property was proved  for exceptional objects
in the category of nilpotent operators of vector spaces with one invariant subspace, where the nilpotency degree is bounded  by $6$
\cite{Dowbor:Meltzer:Schmidmeier} and
   for exceptional objects
in the category of nilpotent operators of vector spaces with two incomparable invariant subspaces,
where the nilpotency degree is bounded  by $3$
\cite{Dowbor:Meltzer:2018}.
Both problems are of tubular type and are related to the Birkhoff problem  \cite{Birkhoff} and to
recent results on stable vector space categories
\cite{Kussin:Lenzing:Meltzer:2013a}, \cite{Kussin:Lenzing:Meltzer:2013b}, \cite{Kussin:Lenzing:Meltzer:2018},

The aim of this paper is to present the following result.
\begin{mainthm}
Let $\La=\La(\pp,\lala)$ be a wild canonical algebra of quiver type, with $\lala=(\la_2,\cdots,\la_t)$. Then \textquotedbl almost all\textquotedbl{} exceptional $\La-$modules can be exhibited by matrices involving as coefficients $\la_i-\la_j$, where $2\leq i,j\leq t$.
\end{mainthm}
The  notion \textquotedbl almost all\textquotedbl \ means that in every $\tauX-$orbit of exceptional modules from a certain place to the right all  modules have the expected matrices.
 We strongly believe that the theorem holds for all exceptional $\La-$modules,
 but the proof of this fact needs  additional arguments.

The theorem will be shown by induction on the rank of a module.
Recall, that matrices for modules of rank $0$ and $1$ are known  \cite{Kussin:Meltzer:2007a}, \cite{Meltzer:2007}.
 Next, by Schofield induction \cite{Schofield}
  each exceptional $\La-$module $M$ of rank greater than or equal to $2$
 can be obtained as the central term of a non-split sequence
$$(\star)\quad 0\lra Y^{\oplus v}\lra M\lra X^{\oplus u}\lra 0,$$
where $(X,Y)$ is an orthogonal exceptional pair in the category $\coh\XX$ of coherent sheaves over the weighted projective line $\XX$
 corresponding to $\La$ and $(u,v)$ is a dimensional vector of an  expceptional representation for the generalized Kronecker algebra having  $\dim_k\ExtX XY$ arrows
\cite{Kedzierski:Meltzer:2013}.
Consequently, like C.M. Ringel in \cite{Ringel:1998} we will study  the category $\Ff(X,Y)$ which consists of all middle term of short exact sequences $(\star)$ for $u,v\in\NN_0$. This category  is equivalent to the module category of generalized Kronecker algebra.
Finally, using an alternative description of  extension spaces we will assign coefficients for exceptional modules over wild canonical algebras.

The result is part of the PhD thesis of the first author at Szczecin University in 2017.
The authors are thankful to C. M. Ringel for helpful discussion concerning the paper \cite{Ringel:1998}.


\section{Notations and basic concepts}

We recall the concept of a weighted projective line in the sense of Geigle-Lenzing \cite{Geigle:Lenzing:1987} associated to a canonical algebra $\La=\La(\pp,\lala)$.
Let $\LL=\LL(\pp)$ be the rank one abelian group with generators
$\vx_1,\dots,\vx_t$
and relations $p_1\vx_1=\cdots p_t\vx_t:=\vc$, where $\vc$ is called \emph{canonical element}.
 Moreover each element $\vy$ of $\LL$ can be written in \emph{normal form} $\vy=a\vc+\sum_{i=1}^t a_i\vx_i$ with $a\in\ZZ$ and $0\leq a_i<p_i$.
The polynomial algebra $k[x_1,\dots,x_t]$ is $\LL-$graded, where degree of $x_i$ is $\vx_i$.
 Because the polynomials $f_i=x^{p_i}_i-x_1^{p_1}-\la_ix_2^{p_2}$ for $i=3,\dots t$ are homogeneous, the quotient algebra $S=k[x_1,\dots,x_t]/\langle f_i\mid i=3,\dots,t\rangle$ is also $\LL-$graded. A \emph{weighted projective line} $\XX$ is by definition the projective spectrum of the $\LL-$graded algebra $S$. The category of coherent sheaves over $\XX$ will be denoted by $\coh\XX$.
In other words the category of coherent sheaves $\coh\XX$ is the Serre quotient
$\mathrm{mod}^{\ZZ}(S)/ \mathrm{mod}^{\ZZ}_0 (S)$, where $\mathrm{mod}^{\ZZ} S$ is the category of finitely generated $\ZZ$-graded modules over $S$ and  $\mathrm{mod}^{\ZZ}_0 (S)  $ the
subcategory of modules of finite length.
It is well known, that each indecomposable sheaf in $\coh\XX$ is a locally free sheaf, called a \emph{vector bundle}, or a \emph{sheaf of finite length}. Denote by $\vect\XX$ (resp. $\cohnull\XX$) the category of vector bundles (resp. finite length sheaves) on $\XX$.

The category $\coh\XX$ is a $\mathrm{Hom}-$finite, abelian $k-$category. Moreover, it is hereditary that means that $\Ext i\XX --=0$ for $i\geq 2$ and it has Serre duality in the form $\ExtX FG\iso D\HomX G{\tauX F}$, where the Auslander-Reiten translation $\tauX$ is given by the shift $F\mapsto F(\vw)$, where $\vw:=(t-2)\vc-\sum_{i=1}^t\vx_i$ denotes the \emph{dualizing element}, equivalently the category
$\coh\XX$ has Auslander-Reiten sequences.
Moreover, there is a tilting object composed of line bundles  
 with $\EndX T=\La$ and it induces an equivalence of a bounded derived category $\mathcal{D}^b(\coh\XX)\stackrel{\iso}{\lra}\mathcal{D}^b(\mod\La)$.

For coherent sheaves there are well known invariants \emph{rank}, \emph{degree} and \emph{determinant}, which correspond to linear forms $\rk, \deg : K_0(\XX)\lra\ZZ$ and $\det:K_0(\XX)\lra \LL(\pp)$, again called rank, degree and determinant.

Recall that a coherent sheaf $E$ over $\XX$ is called \emph{exceptional} if  $\ExtX EE=0$ and $\EndX E$  is a division ring, in case
$k$ is algebraically closed the last  means that $\EndX E=k$.
 A pair $(X,Y)$ in $\coh\XX$ is called \emph{exceptional} if $X$ and $Y$ are exceptional and $\HomX YX=0=\ExtX YX$.
 Finally, an exceptional pair is \emph{orthogonal} if additionally $\HomX XY=0$.

The \emph{rank} of a $\La-$module is defined $\rk M:= \dim_k M_{0}-\dim_k M_{\vc}$.
The rank of a module in this sense equals the rank of the corresponding sheaf in the geometric meaning.
 We denote by $\modplus\La$ (respectively $\modminus\La$ or $\modzero\La$) the full subcategory
 consisting of all $\La-$modules, which indecomposable summands of the decomposition into a
 direct sum have positive (respectively negative or zero) rank. Similarly, by $\cohplus\XX$ (resp. $\cohminus\XX$)
  we denote the full subcategory of all vector bundles over $\XX$, such that the functor $\ExtX T-$ (resp. $\HomX T-$) vanishes.
   Under the equivalence $\mathcal{D}^b(\coh\XX)\stackrel{\iso}{\lra}\mathcal{D}^b(\mod\La)$
\begin{itemize}
\item $\cohplus\XX$ corresponds to $\modplus\La$ by means of $E\mapsto \HomX TE$,
\item $\cohnull\XX$ corresponds to $\modzero\La$ by means of $E\mapsto \HomX TE$,
\item $\cohminus\XX[1]$ corresponds to $\modminus\La$ by means of $E[1]\mapsto \ExtX TE$, where $[1]$ denotes suspension functor of the triangulated category $D^b(\coh\XX)$.
\end{itemize}
For simplicity we will often identify a sheaf $E$ in $\cohplus\XX$ or $\cohnull\XX$ with the corresponded $\La-$module $\HomX TE$.


\section{Exceptional modules of the small rank}\label{sec:samll_rank}
First, we start with some matrix notations. For a natural numbers $n$ by $I_n$ denote the square diagonal matrix of degree $n$ with each non-zero element equal $1$. For a natural number $n$ and $k$ by $X_{n+k,n}$ and $Y_{n+k,n}$ we denote the  following matrices.
$$X_{n+k,n}:=\left[\begin{array}{ccc}
&&\\
&I_n&\\
&&\\
\hline
0&\cdots&0\\
\vdots&\ddots&\vdots\\
0&\cdots&0\\
\end{array}\right]\in M_{n+k,n}(k),\quad
Y_{n+k,n}:=\left[\begin{array}{ccc}
0&\cdots&0\\
\vdots&\ddots&\vdots\\
0&\cdots&0\\
\hline
&&\\
&I_n&\\
&&\\
\end{array}\right]\in M_{n+k,n}(k).$$

A $\La-$module of rank zero is called \emph{regular}. It is well known that the Auslander-Reiten quiver of the regular $\La-$modules consists of a family of orthogonal regular tubes with $t$ exceptional tubes $\Tt_1,\dots, \Tt_t$ of rank $p_1,\dots, p_t$, respectively, while the other tubes are homogeneous. Moreover an exceptional regular modules lies in an exceptional tube and its quasi-length is less of the rank of the tube. We will use the description from \cite{Kussin:Meltzer:2007a}, for the indecomposable regular modules.
However we will only quote the shape of the exceptional ones,
which lies in the tube $\Tt_i$ for $i \in \{3,...,t \} $.
For the tubes $\Tt_1$ and $\Tt_2$ the description is similar.
Following the notations from \cite{Kussin:Meltzer:2007a}
we denote a regular module by $S_a^{[l]}$, where $l$ is the quasi-length of $S_a^{[l]}$ and $a$
indicates the position on the corresponding floor of the tube.
For an exceptional module $S_a^{[l]}$ the quasi-length $l< p_i$ and so all vector space of $S_a^{[l]}$ are zero or one dimensional.

There are $3$ cases:
\begin{itemize}
\item[$(1)$] $1\leq a<p_i\quad\textnormal{and}\quad 0<l<p_i-a$,
\item[$(2)$] $1\leq a<p_i\quad\textnormal{and}\quad p_i-a<l<p_i$,
\item[$(3)$] $a=p_i\quad\textnormal{and}\quad 0<l<p_i$.
\end{itemize}

{\bf Case $(1)$}.

Then $S_a^{[l]}$ has the form
$$\xymatrix @R -15pt {& 0\ar[ldd] & 0\ar[l] && \ar[ll]\cdots& &\ar[ll]0& \ar[l]0 &  \\
& 0\ar[ld] & 0\ar[l] && \ar[ll]\cdots& &\ar[ll]0& \ar[l]0 &\\
0 &&  & & &&&&0, \ar[ld]  \ar[ul]\ar[uul]\ar[ddl]\\
&  \cdots\ar[lu] &\ar[l]0 &\ar[l] k & \cdots\ar[l]_-{\id} &\ar[l]_-{\id} k&\ar[l]0  & \cdots\ar[l] \\
& 0\ar[luu] & 0\ar[l] && \ar[ll]\cdots& &\ar[ll]0& \ar[l]0 &  \\
}$$
where $0\longleftarrow k$ and $k\longleftarrow 0$ correspond to the arrow $\al_a^{(i)}$ and $\al_{a+l}^{(i)}$ in the $i-$th arm.

{\bf Case $(2)$}.

 Let $s:=l-(p_i-a)$. Then $S_a^{[l]}$ is the form
$$
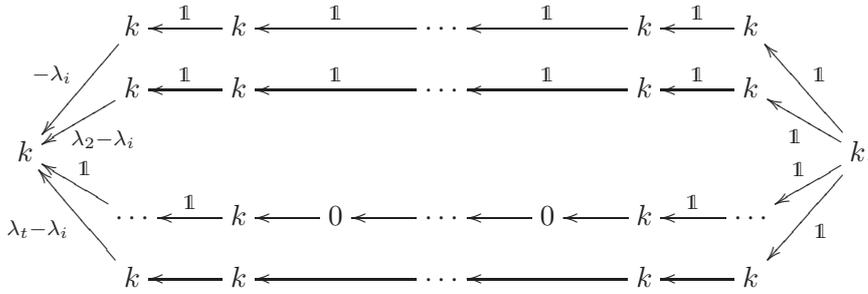
\xymatrix @R -15pt{& k\ar[ldd]_-{-\la_i} & k\ar[l]_-{\id}  && \ar[ll]_-{\id} \cdots& &\ar[ll]_-{\id} k& \ar[l]_-{\id} k &  \\
& k\ar[ld]^-{\la_2-\la_i} & k\ar[l]_-{\id} && \ar[ll]_-{\id} \cdots& &\ar[ll]_-{\id} k& \ar[l]_-{\id} k &\\
k &&  & & &&&&k, \ar[ld]_-{\id}  \ar[ul]^-{\id}\ar[uul]_-{\id}\ar[ddl]^-{\id}\\
&  \cdots\ar[lu]_-{\id} &\ar[l]_-{\id} k &\ar[l] 0 & \cdots\ar[l] &\ar[l]0 &\ar[l]k  & \cdots\ar[l]_-{\id} \\
& k\ar[luu]^-{\la_t-\la_i} & k\ar[l] && \ar[ll]\cdots& &\ar[ll]k& \ar[l]k &  \\
}$$
where $k\longleftarrow 0$ and $0\longleftarrow k$ correspond to the arrow $\al_s^{(i)}$ and  $\al_{a}^{(i)}$ in the $i-$th arm.

{\bf Case $(3)$}
Then $S_a^{[l]}$ is the form
$$
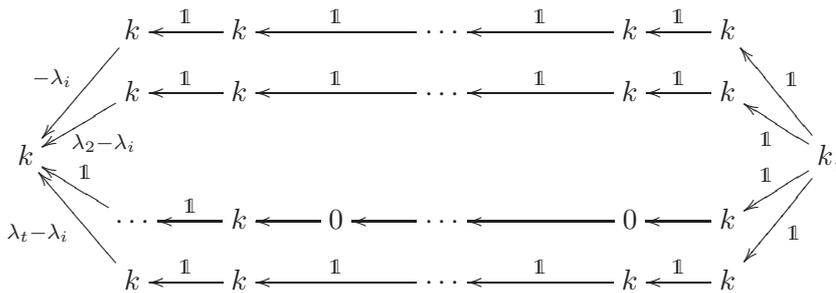
\xymatrix @R -15pt{& k\ar[ldd]_-{-\la_i} & k\ar[l]_-{\id} && \ar[ll]_-{\id} \cdots& &\ar[ll]_-{\id} k& \ar[l]_-{\id} k &  \\
& k\ar[ld]^-{\la_2-\la_i} & k\ar[l]_-{\id} && \ar[ll]_-{\id} \cdots& &\ar[ll]_-{\id} k& \ar[l]_-{\id} k &\\
k &&  & & &&&&k, \ar[ld]_-{\id}  \ar[ul]^-{\id} \ar[uul]_-{\id} \ar[ddl]^-{\id}\\
&  \cdots\ar[lu]_-{\id} &\ar[l]_-{\id} k &\ar[l] 0 & \cdots\ar[l] & &\ar[ll]0  & k\ar[l] \\
& k\ar[luu]^-{\la_t-\la_i} & k\ar[l]_-{\id} && \ar[ll]_-{\id}\cdots& &\ar[ll]_-{\id} k& \ar[l]_-{\id} k &  \\
}$$
where $k\longleftarrow 0$ and $0\longleftarrow k$ correspond to the arrow $\al_l^{(i)}$ and  $\al_{p_i}^{(i)}$ in the $i-$th arm.

For $\La-$modules of rank one there is the following characterization.

\begin{Prop}[\cite{Meltzer:2007}]
\label{thm:case:rank:one}
Let $\La$ be a canonical algebra of  quiver type
and of arbitrary representation type and $M$ an exceptional $\La$-module of rank $1$. Then $M$ is isomorphic to one of the following modules.
$$
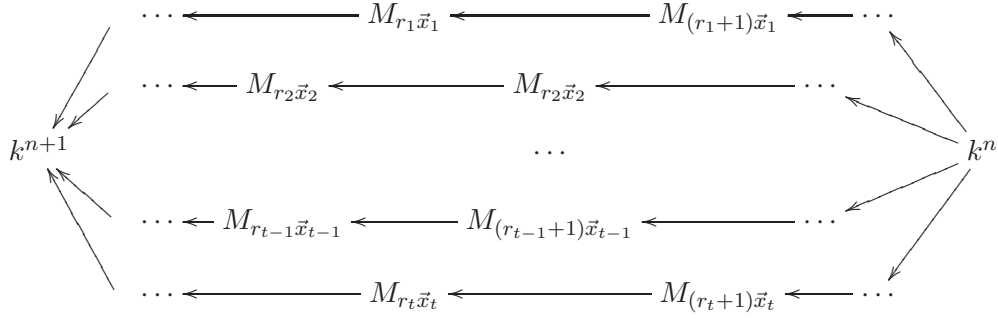
\xymatrix @C -1.8pc @R -15pt{
&&\ar[lldd]&\cdots&&&\ar[lll] M_{r_1\vx_1}  &&  \ar[ll] M_{(r_1+1)\vx_1}  && \ar[ll]\cdots&&\\
&&\ar[lld]&\cdots&&\ar[ll] M_{r_2\vx_2}    & & \ar[ll] M_{r_2\vx_2} & & \ar[ll]\cdots&&\\
k^{n+1}&&& & & & &\cdots & &&&&&\ar[llluu] \ar[llllu] k^n \ar[llldd] \ar[lllld]\\
&&\ar[llu]&\cdots&&\ar[ll] M_{r_{t-1}\vx_{t-1}} & & \ar[ll] M_{(r_{t-1}+1)\vx_{t-1}}  & &\ar[ll]\cdots&&\\
&&\ar[lluu]&\cdots&&&\ar[lll] M_{r_t\vx_t}   & &\ar[ll] M_{(r_t+1)\vx_t} & &\ar[ll]\cdots&&\\
}$$
where $r_i$ is an integer number such that $0\leq r_i<p_i$ for each $i=1,2,\dots,t$ and
\begin{itemize}
\item $M_{s\vx_i}=\left\{\begin{array}{ccl}
k^{n+1} &\text{for} &0\leq s\leq r_i\\
k^{n} &\text{for} &r_i< s\leq p_i\\
\end{array}\right. $
\end{itemize}
Further the matrices of $M$ are given as follows
\begin{itemize}
\item $M_{\al_{s}^{(i)}}=\left\{\begin{array}{ccl}
I_{n+1} &\text{for}& 1<s<r_i\\
I_n &\text{for}& r_i<s\leq p_i
\end{array}\right.\quad\text{for}\quad i=1,2,\dots,t. $
\item $M_{\al_{r_1}^{(1)}}=X_{n+1,n}$
\item $M_{\al_{r_2}^{(2)}}=Y_{n+1,n}$

\end{itemize}
and for $i=3,\dots, t$  we  distinguish
two cases\\
a) $r_i=0$

\begin{itemize}
\item $M_{\al_{1}^{(i)}}=\left[\begin{array}{cccc}
1&0&\cdots& 0\\
\la_i&1&\cdots&0\\
\vdots&\ddots&\ddots&\vdots \\
0&0&\cdots&1\\
0&0&\cdots& \la_i\\
\end{array}\right]\in M_{n+1,n}(k).$
\end{itemize}
b) $r_i>0$
\begin{itemize}
\item $M_{\al_{1}^{(i)}}=\left[\begin{array}{ccccc}
1&0&\cdots&0& 0\\
\la_i&1&\cdots&0&0\\
\vdots&\ddots&\ddots&\vdots &\vdots\\
0&0&\cdots&1&0\\
0&0&\cdots&\la_i&1\\
\end{array}\right]\in M_{n+1}(k)$,\qquad $\bullet$ $M_{\al_{r_i}^{(i)}}=X_{n+1,n}.$
\end{itemize}
\end{Prop}


\section{Schofield induction from sheaves to modules}

Let $M$ be an exceptional object from  $\modplus\La$  of rank greater than or equal to $2$.
Then there is a short exact sequence
$$ (\star)\quad 0\lra Y^{\oplus v} \lra M \lra X^{\oplus u}\lra 0,$$
where $(X,Y)$ is an orthogonal exceptional pair in the category $\coh\XX$, such that the $\rk X< \rk M $ and $\rk Y< \rk M $ and that $(u,v)$
is a dimension vector of an exceptional representation of the generalized Kronecker algebra given by the quiver $\Theta(n)$:
$$\xymatrix {1 \ar @/^1pc/[rr]^{\al_1} \ar @{}[rr]|{\vdots} \ar @/_1pc/[rr]_{\al_n} &&2,}$$
with $n:=\dim_k\ExtX XY$ arrows.

This result is called Schofield induction \cite{Schofield} and was applied
by C. M. Ringel in the situation of exceptional representations
over finite acyclic quivers, hence of hereditary algebras
\cite{Ringel:1998}.

In the case, that the rank of $X$ or $Y$ is at least $2$, we can reapply  Schofield induction again  and as a result we receive the following sequences
$$0\lra Y_2^{\oplus v_2}\lra Y\lra X_2^{\oplus u_2}\lra 0 \quad \text{and}\quad
0\lra Y_3^{\oplus v_3}\lra X\lra X_3^{\oplus u_3}\lra 0.$$

Because with each successive use of the Schofield induction, the rank of the sheaves decreases, after a finite number of steps we receive pairs of exceptional sheaves of rank $0$ or $1$.

 This situation is  illustrated by the following diagram, which has the shape of a tree like the following one.
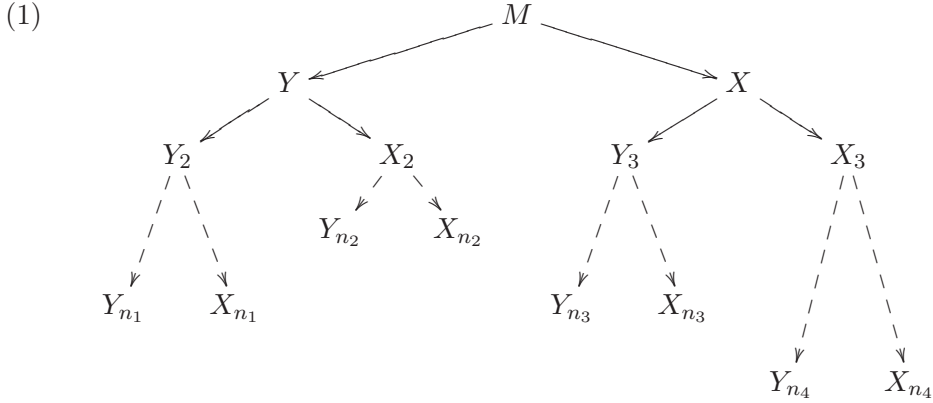
\begin{equation}
	\label{eq:induction_tree}
\xymatrix @C -2pc @R-1pc{
&&&&&&&M\ar @{->}[lllld]\ar @{->}[rrrrd]\\
&&& Y\ar @{->}[lld]\ar @{->}[rrd] &&&&&&&& X\ar @{->}[lld]\ar @{->}[rrd] \\
& Y_2\ar @{-->}[ldd]\ar @{-->}[rdd]& & & & X_2\ar @{-->}[ld]\ar @{-->}[rd]& & & & Y_3\ar @{-->}[ldd]\ar @{-->}[rdd]& & & & X_3\ar @{-->}[lddd]\ar @{-->}[rddd]\\
& && &Y_{n_2}& &X_{n_2}&  & & && & & &\\
Y_{n_1}& &X_{n_1}& && &&  & Y_{n_3}& &X_{n_3}& & & &\\
& && && &&  & & && & Y_{n_4}& &X_{n_4}\\
}
\end{equation}

Applying the functor $\ExtX T-$ to the exact sequence $(\star)$
we see  that if $M$ is a $\La-$module, then each sheaf $X_{i_n}$ such that there is a path form $M$ to $X_{i_n}$ in the tree \eqref{eq:induction_tree} is also $\La-$module. However, we do not know that  a sheave $Y_*$ is a $\La-$module.

The following lemma will allows us, by using the $\tauX$-translation, to shift the tree \eqref{eq:induction_tree} such that all its components will be $\La-$modules.

\begin{Lem}
	\label{lem:translation_line_bundle}
Let $\{L_1,\dots,L_m\}$ be a family of line bundles over $\XX$. Then there is a natural number $N$ such that
$\ExtX T{\tauX^n L_j}=0$ for $j=1,\dots,m$ and for all $n>N$.
\end{Lem}

\begin{pf}
Let $L_j=\Oo\left( a_{j}\vc+\sum_{i=1}^{t}b_{j,i}\vx_i\right)$,
where $a_j\in\ZZ$, $0\leq b_{j,i}\leq p_i-1$ for $j=1,\dots,m$ and $i=1,\dots,t$. We put $N:=\max\Big\{\big\lfloor (1-a_j)(t-2)\big\rfloor+1\mid {1\leq j\leq m}\Big\}$. Then $\vc+\vw-\det \tauX^n L_j=\big(1-a_j-(n-1)(t-2)\big)\vc+\sum_{i=1}^t\big(n-1-b_{j,i} \big)\vx_i<0$
for all $n>N$ and $j=1,\dots,m$.
Therefore by Serre duality
$$(\triangle)\quad\ExtX {\Oo(\vc)}{\tauX^n L_j}\iso D\HomX {\tauX^n L_j}{\Oo(\vc+\vw)}=0 \quad \text{for}\quad j=1,\dots,m.$$
We have to shown that if $n>N$ then
$\ExtX{\Oo(\vx)}{\tauX^n L_j}=0$ for $0\leq\vx <\vc$ and $j=1,\dots,m,$.
Suppose, that $\ExtX{\Oo(\vx)}{\tauX^nL_j}\neq 0$ for some $0\leq\vx<\vc$. Then  using Serre duality we get $\det \tauX^nL_j\leq \vx+\vw.$ Because $\vx+\vw\leq \vc+\vw$, then $\det \tauX^nL_j\leq \vc+\vw$ and
$\ExtX {\Oo(\vc)}{\tauX^n L_j}\iso D\HomX {\tauX^n L_j}{\Oo(\vc+\vw)}\neq 0$
 it is contradictory to $(\triangle)$.
\end{pf}

Immediately from the  lemma above  we receive the following corollary.
\begin{Cor}
	\label{cor:translation_induction_tree}
There is a natural number $N$ such that for $n>N$ all components of the tree \eqref{eq:induction_tree} shifted by $\tauX^n$ are $\La-$modules.
\end{Cor}

\begin{pf}
First,  note that if the sheaves $X$ and $Y$ in the sequence $(\star)$ are $\La-$modules, then middle term $M$ is a $\La-$module. Next, because there
are no nonzero  morphisms from finite length sheaves to vector bundles, each finite length sheaf is a $\La-$module.

Let $\Ll=\{L_i\}_{i\in I}$ be the set of all line bundles appearing in the tree \eqref{eq:induction_tree}. From lemma \ref{lem:translation_line_bundle} applied to the family $\Ll$, there is natural number $N$, such that for all natural number $n>N$ the line bundles $\tauX^nL_i$ are $\La-$modules for $i\in I$. So the vector bundles in the penultimate parts of the tree \eqref{eq:induction_tree} are also $\La-$modules. Moving up from the bottom, we get all sheaves in the $\tauX^n$ image are $\La-$modules.
\end{pf}


\section{Description of extension spaces}

Let $(X,Y)$ be an  orthogonal exceptional  pair in the category $\coh\XX$, this means that
$\HomL XY=0=\HomL YX$, $\ExtL YX=0$ and $\ExtL XY=k^n$ is non zero space.
Assume further that both sheaves $X$ and $Y$ in the sequence $(\star)$ are $\La-$modules.

 We consider the  category $\Ff(X,Y)$, consisting of all right $\La-$modules $M$, that appear as the middle term in a short exact sequence
$$0\lra Y^{\oplus v}\lra M\lra X^{\oplus u}\lra 0\quad\text{for some}\quad v,u\in \mathbb{N}_0.$$
It is well known, that the category $\Ff(X,Y)$ is abelian and has only two simple objects $X$ and $Y$, where the first one  is injective simple and the second one is projective simple \cite{Ringel:1976}.

Acting like C. M. Ringel in the situation of modules over a hereditary algebras \cite{Ringel:1998}
 we show that the problem of classifying the objects in the categories $\Ff(X,Y)$ can be reduced to the classification of the modules over the generalized Kronecker algebra, with $n$ arrows.

To do so let $\eta_1$, \dots ,  $\eta_n$ be  a basis of the vector space  $\ExtL XY$. Thus we have short exact sequences
$$
\eta_i:\quad 0\lra Y\lra Z_i\lra X \lra 0\quad \text{for}\quad i=1,2,\dots n
$$
 From the  ''pull-back'' construction  there is commutative diagram
$$\xymatrix{0\ar[r] & Y^{\oplus n} \ar[r]\ar @{=}[d] & Z \ar[r]\ar[d] & X\ar[r] \ar[d]^{[1_X,\cdots,1_X]^T}& 0\\
0\ar[r] & Y^{\oplus n} \ar[r] & \bigoplus_{i=1}^nZ_i \ar[r] & X^{\oplus n}\ar[r] & 0},$$
where the upper sequence is a universal extension and $Z$ is an exceptional projective object in $\Ff(X,Y)$. In addition, the projective module $Y\oplus Z$ is progenerator of $\Ff(X,Y)$. Therefore the functor $\HomX{Y\oplus Z}-$ induces an equivalence between the category $\Ff(X,Y)$ and the category of modules over the endomorphism algebra $\End\La{Y\oplus Z}$, which is isomorphic to generalized Kronecker algebra $k\Theta(n)$, where $n:=\dim_k\ExtL XY$.

Now, we need a more precise description of the above equivalence.
Recall from \cite{Ringel:1998}  the following concept of extension space between two quiver representations $X$ and $Y$.
Let $C^0(X,Y)$ and $C^1(X,Y)$ be the vector spaces defined as follows
\begin{equation}\nonumber
\begin{split}
C^0(X,Y):=&\bigoplus_{0\leq \vx\leq \vc}\Hom {k}{X_\vx}{Y_\vx}, \\
C^1(X,Y):=&\bigoplus_{\al:\vx \to \vy}\Hom {k}{X_{\vy} }{Y_{\vx}},
\end{split}
\end{equation}
and let $\delta_{X,Y}:C^0(X,Y)\lra C^1(X,Y)$ be the linear map, defined by
$$\delta_{X,Y}\left(\big[f_\vx\big]_{0\leq \vx\leq \vc}\right)=\left[f_{\vy}X_{\al}- Y_{\al} f_{\vz}\right]_{\al:\vy\to \vz,}$$
where  $\al$ passing the set $Q_1$.

For a path algebra $kQ$  the map $\delta_{X,Y}:C^0(X,Y)\lra C^1(X,Y)$ gives also useful description of the extension space of $kQ-$modules \cite{Ringel:1998}. Indeed, then there is $k-$linear isomorphism
$$\Ext 1{kQ}XY\iso C^1(X,Y)/\im (\delta_{X,Y}).$$

For modules over a canonical algebra $\La=\La(\pp,\lala)$ we must additionally consider the canonical relations of the algebra $\La$. For this we take the subspace $U(X,Y)$ of $C^1(X,Y)$ containing all $\left[f_\al\right]_{\al\in Q_1} $ satisfying the following equations.
\begin{equation}\nonumber
	\label{eq:condition_U_X_Y}
\begin{split}
&Y_{\omega_{1,p_i-1}^{(i)}}f_{\al_{p_i}^{(i)}}
+Y_{\omega_{1,p_i-2}^{(i)}}f_{\al_{p_i-1}^{(i)}}X_{\al_{p_i}^{(i)}}
+\cdots + Y_{\al_1^{(i)}}f_{\al_{2}^{(i)}}X_{\omega_{3,p_i}^{(i)}}+ f_{\al_{1}^{(i)}}X_{\omega_{2,p_i}^{(i)}}=\\
=&Y_{\omega_{1,p_1-1}^{(1)}}f_{\al_{p_1}^{(1)}}
+Y_{\omega_{1,p_1-2}^{(1)}}f_{\al_{p_1-1}^{(1)}}X_{\al_{p_1}^{(1)}}
+\cdots + Y_{\al_1^{(1)}}f_{\al_{2}^{(1)}}X_{\omega_{3,p_1}^{(1)}}+ f_{\al_{1}^{(1)}}X_{\omega_{2,p_1}^{(1)}}+\\
+&\lambda_i\left( Y_{\omega_{1,p_2-1}^{(2)}}f_{\al_{p_2}^{(2)}}\right.
+\left. Y_{\omega_{1,p_2-2}^{(2)}}f_{\al_{p_2-1}^{(2)}}X_{\al_{p_2}^{(2)}}
+\cdots + Y_{\al_1^{(2)}}f_{\al_{2}^{(2)}}X_{\omega_{3,p_2}^{(2)}}+ f_{\al_{1}^{(2)}}X_{\omega_{2,p_2}^{(2)}} \right)\\
&\text{for}\quad i=3,4,...,t.
\end{split}
\end{equation}

\begin{Lem}[\cite{Meltzer:2007}]
\label{lem:przedstawienie_Ext}
$\Ext 1{\Lambda}XY\iso U(X,Y)/\im(\delta_{X,Y})$.
\end{Lem}

We recall  the definition of the isomorphism above. Choosing the  bases of the spaces $M_\vx$ we can assume, that for each arrow $\al:i\lra j$ corresponding map $M_{\al}: M_j\lra M_i$ have the shape
$\left[\begin{array}{c|c}
Y_\al&\varphi_\al\\
\hline
0&X_\al\\
\end{array}\right]$. Then an isomorphism $\phi:\ExtL XY\lra U(X,Y)/\im(\delta_{X,Y})$ is given by the formula $M=(M_i,M_\al)\stackrel{\phi}{\mapsto} (\varphi_\al)_{\al\in Q_0}+\im(\delta_{X,Y})$.

Now, we can describe $\La-$modules contained in $\Ff(X,Y)$,
 using the matrices of  $X$, $Y$ and the representation of the quiver $\Theta(n)$, which corresponds to the module $M$. Each module in $\Ff(X,Y)$ can be identified with an element of the extension space $\ExtL{X^{\oplus u}}{Y^{\oplus v}}$. Because $X^{\oplus u}= X\otimes k^u$ and $Y^{\oplus v}= Y\otimes k^v$, then the  space
$\ExtL{X^{\oplus u}}{Y^{\oplus v}}=\ExtL{X\otimes k^u}{Y\otimes k^v}$
is given by the map $\delta_{X\otimes k^u,Y\otimes k^v}$, where the tensor product is taken over the field $k$.
 In this situation the vector space $C^1(X\otimes k^u,Y\otimes k^v)=C^1(X,Y)\otimes\Hom {k}{k^u}{k^v}$ and also $U(X\otimes k^u,Y\otimes k^v)=U(X,Y)\otimes\Hom {k}{k^u}{k^v}$. Therefore, from lemma \ref{lem:przedstawienie_Ext} and from the commutativity of the following diagram
$$\xymatrix{C^0(X\otimes k^u,Y\otimes k^v)\ar @{=}[d]\ar[rr]^-{\delta_{X\otimes k^u,Y\otimes k^v}}&&U(X\otimes k^u,Y\otimes k^v)\ar @{=}[d] \\
C^0(X,Y)\otimes\Hom {k}{k^u}{k^v}\ar[rr]^-{\delta_{X,Y}\otimes 1_{\Hom {}{k^u}{k^v}}} && U(X,Y)\otimes\Hom {k}{k^u}{k^v} }$$
we obtain that
$\ExtL{X\otimes k^u}{Y\otimes k^v}\iso \ExtL XY\otimes \Hom {k}{k^u}{k^v}$.

Let $\phi_1,...,\phi_t$ be a basic of the space $U(X,Y)$.
Then  $\phi_1+\im(\delta_{X,Y}) ,\dots, \phi_n+\im(\delta_{X,Y})$ form a basis of $\ExtL XY$.
Now any element in $\ExtL{X^{\oplus u}}{Y^{\oplus v}}$ is given by an expression $\sum\limits_{k=1}^n\left(f_\al^{(k)}\otimes A_k\right)$, where $A_k\in\Hom k{k^u}{k^v}$ and $\phi_k=\left[f^{(k)}_\al\right]_{\al\in Q_1}$. Therefore an exceptional $\La-$module $M$, that appears in the sequence
$0\lra Y^{\oplus v}\lra M\lra X^{\oplus u}\lra 0,$
has the form
$$M=\left(Y_{\vx}^{\oplus v}\oplus X_{\vx}^{\oplus u},  \left[\begin{array}{c|c}
Y^{\oplus v}_\al&\varphi_\al\\
\hline
0&X^{\oplus u}_\al\\
\end{array}\right] \right)_{0\leq \vx\leq \vc,\  \al\in Q_1}\text{where}\quad\varphi_\al=\sum\limits_{m=1}^n \left(f_{\al}^{(m)}\otimes A_m\right),$$
for an exceptional $\Theta(n)-$representation
 $\xymatrix  {k^v &&\ar @/^1pc/[ll]^{A_n} \ar @{}[ll]|{\vdots} \ar @/_1pc/[ll]_{A_1} k^u}$ of the generalized Kronecker algebra.
An explicit  basis for the subspace $U(X,Y)$ we will construct in the next section.

 Now we will focus on exceptional modules over the generalized Kronecker algebra.
 The exceptional modules in this case are known. They are preprojective or preinjective and can be exhibited 
 by matrices having only $0$ and $1$ entries \cite{Ringel:1998} 
For recent results concerning modules over generalized Kronecker algebra we refer to
\cite{Ringel:2013}, \cite{Ringel:2018}, \cite{Ringel:2016} \cite{Weist}.

\begin{Lem}\label{lem:Kronecker_niezerowe_wyrazy}
Let  $V=\xymatrix {k^v &&\ar @/^1pc/[ll]^{A_n} \ar @{}[ll]|{\vdots} \ar @/_1pc/[ll]_{A_1} k^u}$ be an exceptional representation of the quiver $\Theta(n)$ and let $A_m=\left[a_{i,j}^{(m)}\right]$ for $m=1,\dots,n$.
Then for each pair of natural numbers $(i,j)$ there is at most one index $m$ such that the coefficient $a_{i,j}^{(m)}$ of the matrix $A_m$ is non-zero.
\end{Lem}

\begin{pf} We will use the description of the extension space to show that if for two matrices $A_1$ and $A_2$ of an exceptional
 representation $V$  of $\Theta(n)$
 non-zero coefficient  appear at the same row and column, then $\Ext 1{k\Theta(n)}VV\neq 0$.
Consider the map
$\delta=\delta_{V,V}:C^0(V,V)\lra C^1(V,V)$, where
$C^0(V,V)=\Hom k{k^v}{k^v}\bigoplus \Hom k {k^u}{k^u}$ and  $C^1(V,V)=\bigoplus_{m=1}^n\Hom k{k^u}{k^v}$.

Then for $(f,g) \in C^0(V,V)$ we have
$\delta(f,g)=\bigoplus_{m=1}^n(f A_m-A_m g)$.
The vector space  $C^0(V,V)$ has a base of the form
$(e_{i,j}^v,0)$ for $1\leq i,j\leq v$, and  $(0,e_{i,j}^u)$ for $1\leq i,j\leq u,$
where  $e_{i,j}^*$ is an elementary matrix with one non-zero element (equal $1$) in the $i\times j-$place. Because $A_m=\left[a_{i,j}^{(m)}\right]$ for $m=1,\dots,n$, then
 $\im(\delta)$ is generated by the elements
\begin{equation}\nonumber
\begin{split}
\delta(e_{i,j}^v,0)&=\bigoplus_{m=1}^n e_{i,j}^vA_m=\bigoplus_{m=1}^n\left[\begin{array}{ccccccc}
0&\cdots&0&a_{1,j}^{(m)}&0&\cdots&0\\
\vdots&&\vdots&\vdots&\vdots&&\vdots\\
0&\cdots&0&a_{v,j}^{(m)}&0&\cdots&0\\
\end{array}\right]\\
\delta(0,e_{i,j}^u)&=\bigoplus_{m=1}^n A_me_{i,j}^u=\bigoplus_{m=1}^n\left[\begin{array}{cccc}
0&0&\cdots &0\\
\vdots&\vdots& &\vdots\\
0&0&\cdots &0\\
a_{i,1}^{(m)}&a_{i,2}^{(m)}&\cdots&a_{i,u}^{(m)}\\
0&0&\cdots &0\\
\vdots&\vdots& &\vdots\\
0&0&\cdots &0\\
\end{array}\right].
\end{split}
\end{equation}
Without lost of generality we can assume that $a_{1,1}^{(m)}\neq 0$  for $m=1,2$. Then the element $x=e_{1,1}\oplus 0\oplus\cdots\oplus 0$ belongs to $C^1(V,V)$ and   $x\notin\im(\delta)$. Therefore $\Ext 1{k\Theta(n)}VV\iso C^1(V,V)/\im(\delta)\neq 0$.
\end{pf}


\section{A construction of a base for $U(X,Y)$.}

Let $\La$ be a canonical algebra of the type $\pp=(p_1,\dots,p_t)$
and with parameters $\lala=(\la_2=0,\la_3=1,\dots, \la_t)$. A representation $M=\rep MQ{\vx}{\al}$ of an exceptional
 $\La-$module with positive rank is called \emph{acceptable} if it satisfies the following conditions.
\begin{itemize}
\item[\textbf{C1.}] The matrices $M_{\al_{1}^{(1)}}$, $M_{\al_{1}^{(3)}}$, $M_{\al_{1}^{(4)}}$, \dots, $M_{\al_{1}^{(t)}}$ have entries of the form $\lambda_a-\lambda_b$ for same $a,b\geq 2$, only.

\item[\textbf{C2.}] All other matrices have only $0$ and $1$ as their coefficients.

\item[\textbf{C3.}] For each path $\omega_{u,v}^{(2)}:(u-1)\vx_2\lra v\vx_2$  the entries of the matrix $M_{\omega_{u,v}^{(2)}}$ are equal to $0$ or $1$.

\item[\textbf{C4.}] For each path $\omega_{1,v}^{(i)}:0\lra v\vx_i$, where $i\neq 2$
  the entries of the matrix  $M_{\omega_{1,v}^{(i)}}$ are of the form $\lambda_a-\lambda_b$ for same $a,b\geq 2$.

\item[\textbf{C5.}] For each path $\omega_{u,v}^{(i)}:(u-1)\vx_i\lra v\vx_i$, where $i\neq 2$ and $u\geq 2$
the entries of the matrix $M_{\omega_{u,v}^{(i)}}$ are equal to $0$ or $1$.
\end{itemize}


The following lemma \cite[Lemma 3.4]{Meltzer:2007} is useful.
\begin{Lem}
	\label{lem:odwzorania_rkM>0}
Let $M$ be an acceptable representation of a module in $\modplus\La$.
Then by base change we can assume that
$$M_{\al_{j}^{(i)}}=\left[\begin{array}{ccc}
\\
&I_n&\\
\\
\hline
0&\cdots& 0\\
\vdots&\ddots&\vdots\\
0&\cdots &0\\
\end{array}\right]\quad \text{for}\quad 2\leq j\leq p_i,\quad i=3,...,t.$$
In addition, the matrices $M_{\al_1^{(i)}}$ have again
entries $\lambda_a-\lambda_b$ for same $a,b\geq 2$.\kwadracik
\end{Lem}

For an exceptional pair $(X,Y)$ with acceptable representations of $X$ and $Y$
we will construct a basis of subspace $U(X,Y)$,
for which each basis vector has only coefficients of the form $\lambda_a-\lambda_b$.
In the case that the ranks of $\rk X>0$ and $\rk Y>0$ this was done in \cite{Meltzer:2007}.

\begin{Lem}\label{lem:base_for_positive_rank}
 Let $X$ and $Y$ be $\La$-modules in $\modplus\La$, with acceptable representations. Then there is basis $F^{(1)},\dots, F^{(d)}$ of the subspace $U(X,Y),$ where $F^{(j)}=\left[f_\al^{(j)}\right]_{\al\in Q_1}$ satisfies the following properties:
\begin{itemize}
\item[$(i)$]  The entries of the matrix  $f_{\al_j^{(2)}}$ are equal $0$ or $1$ for $1 \leq j\leq p_i$,
\item[$(ii)$] The entries of the matrix   $f_{\al_j^{(i)}}$ are equal $0$, $1$  for  $2 \leq j\leq p_i$ and $i=1,3,4,\dots,t$.
\item[$(iii)$] The entries of the matrix
$f_{\al_1^{(i)}}$ are equal $\la_a-\la_b$ for $i=1,3,4,\dots,t$.
\end{itemize}
\end{Lem}

Note, that in the sequence $(\star)$ of the Schofield induction the $\La-$module $Y$ always have the positive rank, but $X$ can have rank zero. In this situation, we need one more lemma.

\begin{Lem} \label{lem:baze_for_regular}
Let $Y$ and $X$ be an exceptional $\La-$modules such that $Y\in\modplus\La$ and
  $X \in\modzero\La$.
  Assume that $X$ and $Y$ have acceptable representations
  and $X$ lies in the exceptional tube corresponding to $i-$th arm of the canonical algebra.
   Then there is a base $F^{(1)},\dots,F^{(d)}$ of the subspace $U(X,Y)$, where $F^{(j)}=\big[f_\al^{(j)}\big]_{\al\in Q_1}$ satisfies the following properties:
\begin{itemize}
\item[$(i)$] The entries of the matrix
$f_{\al_j^{(2)}}$ are equal $0$, $1$ for $1\leq j\leq p_2$,
\item[$(ii)$] The entries of the matrix
$f_{\al_j^{(i)}}$ are equal $0$, $1$ for $1\leq j\leq p_i$,
\item[$(iii)$] The entries of the matrix
$f_{\al_1^{(m)}}$ are equal $\lambda_a-\lambda_b$ for $m=1,3,4,\dots, t$ $\wedge$ $m\neq i$,
\item[$(iv)$] The entries of the matrix   $f_{\al_j^{(m)}}$ are equal $0$, $1$ for $2\leq j\leq p_m$ and  $m=3,4,\dots, t$.
\end{itemize}
\end{Lem}

\begin{pf} Because $X$ lies in the exceptional tube corresponding to $i-$th arm of the canonical algebra,
 then it has  a representation of the form $S_a^{[l]}$ from the section \ref{sec:samll_rank}, such that
\begin{itemize}
\item[$(1)$] $1\leq a<p_i\quad\textnormal{and}\quad 0<l<p_i-a$,
\item[$(2)$] $1\leq a<p_i\quad\textnormal{and}\quad p_i-a<l<p_i$,
\item[$(3)$] $a=p_i\quad\textnormal{and}\quad 0<l<p_i$.
\end{itemize}
In particular, all vector space of $S_a^{[l]}$ are zero or one dimensional.

{\bf  Case $(1)$}.
From the shape of $S_a^{[l]}$ any element of the subspace $U(X,Y)$ has the form
$$F=\left[\begin{array}{ccccccccc}
0&\cdots& 0&0&\cdots&0&0&\cdots&0\\
\vdots& &\vdots&\vdots&&\vdots&\vdots&&\vdots\\
0&\cdots& 0&0&\cdots&0&0&\cdots&0\\
0& \cdots&0& f_{\al_{a}^{(i)}} & \cdots & f_{\al_{a+l-1}^{(i)}}& 0&\cdots &0 \\
0&\cdots& 0&0&\cdots&0&0&\cdots&0\\
\vdots& &\vdots&\vdots&&\vdots&\vdots&&\vdots\\
0&\cdots& 0&0&\cdots&0&0&\cdots&0\\
\end{array}\right].$$
In addition, the condition describing the subspace $U(X,Y)$ vanishes.

Now we fix $j$ such that $a\leq j < a+l$. Let denote by $e_{r}$  the matrix unit (the matrix with one coefficient $1$ namely the coefficient in the row with index $r$, the remaining coefficients are zero). Then $e_{r}$ is an element in
 $\Hom{k}{X_{j\vec{x}_i}}{Y_{(j-1)\vec{x}_i}}= \Hom{k}{k}{Y_{(j-1)\vec{x}_i}}$, where $1\leq r\leq \dim _kY_{(j-1)\vec{x}_i}$ and
$$F_j^r=\left[\begin{array}{ccccccccc}
0&\cdots& 0&0&0&\cdots&0\\
\vdots& &\vdots&\vdots&\vdots&&\vdots\\
0&\cdots& 0&0&0&\cdots&0\\
0& \cdots&0& e_r & 0&\cdots &0 \\
0&\cdots& 0&0&0&\cdots&0\\
\vdots& &\vdots&\vdots&\vdots&&\vdots\\
0&\cdots& 0&0&0&\cdots&0\\
\end{array}\right].$$
belongs to $U(X,Y)$ ($e_{r}$ lies in $j$-th column). It is easy to check, that $F_j^r$ for $1\leq r\leq \dim _kY_{(j-1)\vec{x}_i}$ and  $a\leq j < a+l$ create a base of the subspace $U(X,Y)$.

{\bf  Case $(1)$}.
Any element of the subspace $U(X,Y)$ has the form
$$F_j^r=\left[\begin{array}{ccccccccc}
f_{\al_{1}^{(1)}}&\cdots& f_{\al_{s-1}^{(1)}}& f_{\al_{s}^{(1)}}&\cdots & f_{\al_{a-1}^{(1)}}&f_{\al_{a}^{(1)}} &\cdots&f_{\al_{p_1}^{(1)}}\\
\vdots& &\vdots&\vdots&\vdots&&\vdots\\
f_{\al_{1}^{(i-1)}}&\cdots& f_{\al_{s-1}^{(i-1)}}& f_{\al_{s}^{(i-1)}}&\cdots & f_{\al_{a-1}^{(i-1)}}&f_{\al_{a}^{(i-1)}} &\cdots&f_{\a l_{p_{i-1}}^{(1)}}\\
f_{\al_{1}^{(i)}}& \cdots& f_{\al_{s-1}^{(i)}}& 0 & \cdots& 0 &f_{\al_{a}^{(i)}} &\cdots&f_{\al_{p_i}^{(i)}} \\
f_{\al_{1}^{(i+1)}}&\cdots& f_{\al_{s-1}^{(i+1)}}& f_{\al_{s}^{(i+1)}}&\cdots & f_{\al_{a-1}^{(i+1)}}&f_{\al_{a}^{(i+1)}} &\cdots &f_{\al_{p_{i+1}}^{(1)}}\\
\vdots& &\vdots&\vdots&\vdots&&\vdots\\
f_{\al_{1}^{(t)}}&\cdots& f_{\al_{s-1}^{(t)}}& f_{\al_{s}^{(t)}}&\cdots & f_{\al_{a-1}^{(t)}}&f_{\al_{a}^{(t)}} &\cdots&f_{\al_{p_t}^{(t)}}\\
\end{array}\right],$$
where the condition described $U(X,Y)$ has the following shape.
\begin{equation}\nonumber
\begin{split}
&Y_{\omega_{1,p_i-1}^{(i)}}f_{\al_{p_i}^{(i)}}+\cdots
+Y_{\omega_{1,a-1 }^{(i)}}f_{\al_{a}^{(i)}}\\
=&Y_{\omega_{1,p_1-1}^{(1)}}f_{\al_{p_1}^{(1)}}
+Y_{\omega_{1,p_1-2}^{(1)}}f_{\al_{p_1-1}^{(1)}}
+\cdots + Y_{\al_1^{(1)}}f_{\al_{2}^{(1)}}+ f_{\al_{1}^{(1)}}\\
+&\lambda_i\left( Y_{\omega_{1,p_2-1}^{(2)}}f_{\al_{p_2}^{(2)}}\right.
+\left. Y_{\omega_{1,p_2-2}^{(2)}}f_{\al_{p_2-1}^{(2)}}
+\cdots + Y_{\al_1^{(2)}}f_{\al_{2}^{(2)}}+ f_{\al_{1}^{(2)}}\right)
\end{split}
\end{equation}
and for $j\in\{3,4,\dots,t\}$ and $j\neq i$ we get
\begin{equation}\nonumber
\begin{split}
&Y_{\omega_{1,p_j-1}^{(j)}}f_{\al_{p_j}^{(j)}}
+Y_{\omega_{1,p_j-2}^{(j)}}f_{\al_{p_j-1}^{(j)}}
+\cdots + Y_{\al_1^{(j)}}f_{\al_{2}^{(j)}}+ f_{\al_{1}^{(j)}}\\
=&Y_{\omega_{1,p_1-1}^{(1)}}f_{\al_{p_1}^{(1)}}
+Y_{\omega_{1,p_1-2}^{(1)}}f_{\al_{p_1-1}^{(1)}}
+\cdots + Y_{\al_1^{(1)}}f_{\al_{2}^{(1)}}+ f_{\al_{1}^{(1)}}\\
+&\lambda_i\left( Y_{\omega_{1,p_2-1}^{(2)}}f_{\al_{p_2}^{(2)}}\right.
+\left. Y_{\omega_{1,p_2-2}^{(2)}}f_{\al_{p_2-1}^{(2)}}
+\cdots + Y_{\al_1^{(2)}}f_{\al_{2}^{(2)}}+ f_{\al_{1}^{(2)}}\right).
\end{split}
\end{equation}
We fix $j$ such that $2\leq j\leq p_i$. Again $e_{r}$ is unit matrix belongs to
$\Hom{k}{X_{j\vx_1}}{Y_{(j-1)\vx_1}}= \Hom{k}{k}{Y_{(j-1)\vx_1}}$ for $1\leq r\leq \dim _kY_{(j-1)\vx_1}$. Then the element
$$F_{\al_j^{(1)}}^{r}=\left[\begin{array}{cccccccc}
-Y_{\omega_{1,j-1}^{(1)}} e_{r} &0& \cdots & 0 & e_{r} &0 &\cdots&0\\
0 &0& \cdots & 0 & 0 &0 &\cdots&0\\
0 & 0& \cdots & 0 & 0 &0 &\cdots&0\\
\vdots & \vdots& \ddots & \vdots & \vdots &\vdots &\ddots&\vdots\\
0 & 0& \cdots & 0 & 0 &0 &\cdots&0\\
\end{array}\right]$$
(where $e_{r}$ lies in the $j$-th column) belongs to $U(X,Y)$.

We fix $j$ such that $1\leq j\leq p_2$ and let $e_{r}$ belongs to $\Hom{k}{X_{j\vec{x}_2}}{Y_{(j-1)\vec{x}_2}}=\Hom{k}{k}{Y_{(j-1)\vec{x}_2}}$ for $1\leq r\leq \dim _kY_{(j-1)\vec{x}_2}$. Then the element
$$F_{\al_j^{(2)}}^{r}=\left[\begin{array}{cccccccc}
-\lambda_i Y_{\omega_{1,j-1}^{(2)}} e_{r} & 0&\cdots  &0 &0&0&\cdots&0\\
0 & 0& \cdots & 0 & e_{r} &0 &\cdots&0\\
(\la_3-\la_i) Y_{\omega_{1,j-1}^{(2)}} e_{r} & 0& \cdots & 0 & 0 &0 &\cdots&0\\
\vdots & \vdots& \ddots & \vdots & \vdots &\vdots &\ddots&\vdots\\
(\lambda_{i-1}-\lambda_i) Y_{\omega_{1,j-1}^{(2)}} e_{r}  & 0& \cdots & 0 & 0 &0 &\cdots&0\\
0 & 0& \cdots & 0 & 0 &0 &\cdots&0\\
(\lambda_{i+1}-\lambda_i) Y_{\omega_{1,j-1}^{(2)}} e_{r}  & 0& \cdots & 0 & 0 &0 &\cdots&0\\
\vdots & \vdots& \ddots & \vdots & \vdots &\vdots &\ddots&\vdots\\
(\lambda_t-\lambda_i) Y_{\omega_{1,j-1}^{(2)}} e_{r}  & 0& \cdots & 0 & 0 &0 &\cdots&0\\
\end{array}\right]$$
belongs to $U(X,Y)$.

Next, assume that $3\leq m\leq t$, $m\neq i$ and $1<j\leq p_m$. Let $e_{r}$ be a unit matrix in $\Hom{k}{X_{j\vec{x}_m}}{Y_{(j-1)\vec{x}_m}}=\Hom{k}{k}{Y_{(j-1)\vec{x}_m}}$ for $1\leq r\leq \dim _kY_{(j-1)\vec{x}_k}$. Then the element
$$F_{\al_j^{(m)}}^{r}=\left[\begin{array}{cccccccc}
0 & 0& \cdots & 0 & 0 &0 &\cdots&0\\
0 & 0& \cdots & 0 & 0 &0 &\cdots&0\\
0 &0& \cdots & 0 & 0 &0 &\cdots&0\\
\vdots & \vdots& \ddots & \vdots & \vdots &\vdots &\ddots&\vdots\\
-Y_{\omega_{1,j-1}^{(m)}}e_{r} & 0& \cdots & 0 & e_{r} &0 &\cdots&0\\
\vdots & \vdots& \ddots & \vdots & \vdots &\vdots &\ddots&\vdots\\
0 & 0& \cdots & 0 & 0 &0 &\cdots&0\\
\end{array}\right]$$
belongs to $U(X,Y)$.

Now, we fix $j$ such that $j\in\{a, a+1,\dots, p_i  \}$ and let $e_{r}$ be a unit matrix in $\Hom{k}{X_{j\vec{x}_i}}{Y_{(j-1)\vec{x}_i}}=\Hom{k}{k}{Y_{(j-1)\vec{x}_i}}$, where $1\leq r\leq \dim _kY_{(j-1)\vec{x}_i}$. Then the element
$$F_{\al_j^{(i)}}^{r}=\left[\begin{array}{cccccccc}
Y_{\omega_{1, j-1}^{(i)}}e_{r} & 0& \cdots & 0 & 0 &0 &\cdots&0\\
0 & 0& \cdots & 0 & 0 &0 &\cdots&0\\
Y_{\omega_{1, j-1}^{(i)}}e_{r}& 0& \cdots & 0 & 0 &0 &\cdots&0\\
\vdots & \vdots& \ddots & \vdots & \vdots &\vdots &\ddots&\vdots\\
Y_{\omega_{1, j-1}^{(i)}}e_{r}& 0& \cdots & 0 & 0 &0 &\cdots&0\\
0 & \cdots& \cdots & 0 & e_{r} &0 &\cdots&0\\
Y_{\omega_{1, j-1}^{(i)}}e_{r} & 0& \cdots & 0 & 0 &0 &\cdots&0\\
\vdots & \vdots& \ddots & \vdots & \vdots &\vdots &\ddots&\vdots\\
Y_{\omega_{1, j-1}^{(i)}}e_{r} & 0& \cdots & 0 & 0 &0 &\cdots&0\\
\end{array}\right]$$
belongs to $U(X,Y)$.

Let $j$ be a natural number such that $j\in\{1, 2,\dots, s-1  \}$ and let $e_{r}$ belongs to $\Hom{k}{X_{j\vec{x}_i}}{Y_{(j-1)\vec{x}_i}}=\Hom{k}{k}{Y_{(j-1)\vec{x}_i}}$, for   $1\leq r\leq \dim _kY_{(j-1)\vec{x}_i}$. Then the element
$$F_{\al_j^{(i)}}^r=\left[\begin{array}{ccccccccc}
0&\cdots& 0&0&0&\cdots&0\\
\vdots& &\vdots&\vdots&\vdots&&\vdots\\
0&\cdots& 0&0&0&\cdots&0\\
0& \cdots&0& e_r & 0&\cdots &0 \\
0&\cdots& 0&0&0&\cdots&0\\
\vdots& &\vdots&\vdots&\vdots&&\vdots\\
0&\cdots& 0&0&0&\cdots&0\\
\end{array}\right],$$
belongs to $U(X,Y)$, where $e_r$ lies in $j-$th column and $i-$th row.

It is easy to check, that $F_{\al}^{r}$ are a base of $U(X,Y)$. In the end, we must check that matrices $f_{\al}^{r}$ of basis vectors $F_{\al}^{r}$ have desired entries. Because the representation for $Y$ is acceptable, then
the matrices $Y_{\omega_{u,v}^{(m)}}$ have only entries $0$, $\pm 1$, $\pm\la_a$, $\la_a-\la_b$. Hence the matrix $Y_{\omega_{1,j-1}^{(m)}}e_{r}$ has the same entries. In addition, for $m=2$, the  coefficients of matrices $Y_{\omega_{u,v}^{(2)}}$ are equal to $0$ or $1$. Therefore matrices  $(\la_a-\la_b)Y_{\omega_{1,j-1}^{(2)}}e_r$ have only entries $0$, $\pm 1$, $\pm \la_a$, $\lambda_a-\lambda_b$.

The case $(3)$ is similar to $(2)$.
\end{pf}

Remark, that the coefficients of the form $\la_a-\la_b$ occur as coefficients of basis vector of $U(X,Y)$ only if they appear in the acceptable representations $X$ or $Y$. In particular if $X$ and $Y$ are rank modules of rank $1$,
 from Proposition \ref{thm:case:rank:one},
  then the all basis vectors of $U(X,Y)$ have only coefficients $0$, $\pm 1$ and $\pm\la_i$.


\section{Proof of the main theorem}

\begin{Prop}[Induction step]
	\label{Thm:induction_step}
Let $M$ be an exceptional module over a canonical algebra $\La$, such that $\rk M\geq 2$. Let $(X,Y)$ be an orthogonal exceptional pair of $\La-$modules, obtained from Schofield induction applied to $M$. If $X$ and $Y$ allows acceptable representations, then also $M$ allows an acceptable representation.
\end{Prop}

\begin{pf}
We will use the basis $F^{(1)}=\left[f_{\al}^{(1)}\right]_{\al\in Q_1}$,\dots, $F^{(n)}=\left[f_{\al}^{(n)}\right]_{\al\in Q_1}$ of the subspace $U(X,Y)$ from the lemma \ref{lem:base_for_positive_rank} or lemma \ref{lem:baze_for_regular}.
Because $M$ belong to $\Ff(X,Y)$, it has the following form.
$$M=\left(Y_{\vx}^{\oplus v}\oplus X_{\vx}^{\oplus u},  \left[\begin{array}{c|c}
Y^{\oplus v}_\al&\varphi_\al\\
\hline
0&X^{\oplus u}_\al\\
\end{array}\right] \right)_{0\leq \vx\leq \vc,\  \al\in Q_1}\text{where}\quad\varphi_\al=\sum\limits_{m=1}^n \left(f_{\al}^{(m)}\otimes A_m\right),$$
for an exceptional $\Theta(n)-$representation
 $\xymatrix  {k^v &&\ar @/^1pc/[ll]^{A_n} \ar @{}[ll]|{\vdots} \ar @/_1pc/[ll]_{A_1} k^u}$.
Recall that all matrices $A_1$,\dots, $A_n$ have entries only $0$ and $1$ ( see \cite{Ringel:1998}) and moreover from lemma \ref{lem:Kronecker_niezerowe_wyrazy} non-zero coefficients in consecutive matrices $A_1$,\dots,$A_n$ occur in different places. Therefore the matrix
$\sum\limits_{m=1}^{n}(f^{(m)}_{\al}\otimes A_m)$ has the same entries as matrices of basis vector $F^{(1)}$,\dots, $F^{(n)}$ of $U(X,Y)$.
Therefore, because  $X$ and $Y$ are acceptable, the matrix $M_{\al_1^{(i)}}$ has entries of the form $\lambda_a-\lambda_b$ for $i=1,3,4,\dots,t$ and for $i=2$ only $0$ and $1$ appear. Next, $M_{\al_j^{(i)}}$  is a zero-one matrix for $2\leq j\leq p_i$ and $i=1,2,\dots,t$.

Now, we must check, that for each path $\omega_{l,m}^{(i)}=\al_m^{(i)}\dots\al_l^{(i)}$ the matrix $M_{\omega_{l,m}^{(i)}}=M_{\al_l^{(i)}}\dots M_{\al_m^{(i)}}$ has only expected coefficients. After standard calculations we obtain, that
$$M_{\omega_{l,m}^{(i)}}= \left[\begin{array}{c|c}
\left(Y_{\omega_{l,m}^{(i)}}\right)^{\oplus v}&\sum YfX_{\omega_{l,m}^{(i)}}\\
\\
\hline
\\
0& \left(X_{\omega_{l,m}^{(i)}}\right)^{\oplus u}
\end{array}\right],$$
where by $\sum YfX_{\omega_{l,m}^{(i)}}$ we denote
$$\sum\limits_{j=1}^{n} \left\{Y_{\omega_{l,m-1}^{(i)}} f^{(j)}_{{\al_m^{(i)}}}+ Y_{\omega_{l,m-2}^{(i)}} f^{(j)}_{\al_{m-1}^{(i)}} X_{\al_{m}^{(i)}}+ ...+
f^{(j)}_{\al_l^{(i)}} X_{\omega_{l+1,m}^{(i)}}\right\}\otimes A_j.$$
Because $X$ and $Y$ are acceptable, then  $Y_{\omega_{l,m}^{(i)}}$ and $X_{\omega_{l,m}^{(i)}}$ allow only desired entries. Again $\sum YfX_{\omega_{l,m}^{(i)}}$ has the same coefficients as $$Y_{\omega_{l,m-1}^{(i)}} f^{(j)}_{{\al_m^{(i)}}}+ Y_{\omega_{l,m-2}^{(i)}} f^{(j)}_{\al_{m-1}^{(i)}} X_{\al_{m}^{(i)}}+ ...+
f^{(j)}_{\al_l^{(i)}} X_{\omega_{l+1,m}^{(i)}}.$$
Now the statement concerning the coefficients of the
matrices $M_{\omega_{l,m}^{(i)}}$ follows from the explicit description of the elements $F^{(i)}$ by a case by case inspection.
\end{pf}

Let us note, that coefficients of the form $\la_a-\la_b$ appear only for regular modules.
This means that if in the tree \eqref{eq:induction_tree} there are only vector bundles, then each modules in this  tree (after translations) can by  established by matrices having coefficients $0$, $\pm 1$, $\pm\la_a$.

\begin{pf}[Proof of Main Theorem]
We prove the fact by induction on the rank of the exceptional module. Remember that a description of  exceptional modules of the zero and one rank in section \ref{sec:samll_rank}, which gives us the start of induction.
Let $M$ be an exceptional $\La-$module of rank r and assume that $r\geq 2$. Then $M$ is corresponds to an exceptional vector bundle over the weighted projective line $\XX$ associated to $\La$. By repeated use of Schofield induction, we obtain the figure \eqref{eq:induction_tree} in the category $\coh\XX$ for $M$. Then from Corollary \ref{cor:translation_induction_tree} we can shift all tree, such that each sheaf in the tree is a $\La-$module. Therefore, up to \textquotedbl almost all\textquotedbl{} we can assume that all tree \eqref{eq:induction_tree} belongs to the category $\mod\La$. Because all tree components have smaller rank than $M$, then they have acceptable representations. Therefore the claim follows from Proposition \ref{Thm:induction_step}.
\end{pf}

\bibliographystyle{amsplain}


\end{document}